\documentclass[11pt]{amsart}

\usepackage[usenames,dvipsnames,svgnames,table]{xcolor}
\usepackage[hyphens]{url}
\usepackage[pagebackref,linktocpage=true,colorlinks=true,
   linkcolor=Blue,citecolor=BrickRed,urlcolor=RoyalBlue]{hyperref}
\usepackage[msc-links,abbrev]{amsrefs}
\usepackage{amsmath,amsthm,amssymb}
\usepackage{mathtools}
\usepackage{enumitem}
\usepackage{overpic}

\newcommand{\C}{\mathcal{C}}
\newcommand{\R}{\mathbf{R}}
\newcommand{\Sph}{\mathbf{S}}
\newcommand{\p}{\partial}
\DeclareMathOperator{\Hess}{Hess}

\newtheorem{theorem}{Theorem}[section]
\newtheorem{lemma}[theorem]{Lemma}
\newtheorem{proposition}[theorem]{Proposition}

\theoremstyle{definition}

\renewcommand{\epsilon}{\varepsilon}

\begin{document}

\title{Nonrigidity of smooth tight surfaces}

\author[M. Ghomi]{Mohammad Ghomi}
\address{School of Mathematics, Georgia Institute of Technology,
  Atlanta, GA 30332}
\email{ghomi@math.gatech.edu}
\urladdr{ghomi.math.gatech.edu}

\author[H. Qiu]{Hongda Qiu}
\address{School of Mathematics, Georgia Institute of Technology,
  Atlanta, GA 30332}
\email{qhongda3@gatech.edu}

\author[M. Raffaelli]{Matteo Raffaelli}
\address{School of Mathematics, Georgia Institute of Technology,
  Atlanta, GA 30332}
\email{raffaelli@math.gatech.edu}
\urladdr{matteoraffaelli.com}

\author[Y. Zhang]{Yuetong Zhang}
\address{School of Mathematical Sciences, Peking University,
  Beijing 100871, China}
\email{2300010810@stu.pku.edu.cn}

\date{\today (Last Typeset)}

\begin{abstract}
We construct smooth isometric embeddings of closed
surfaces of positive genus in Euclidean $3$-space which are tight,
i.e., have minimal total positive curvature, but are not congruent.
These embeddings agree outside a compact subset of the region of
negative curvature, and enclose different volumes. Thus the rigidity
theorem of Alexandrov for analytic tight surfaces does not hold in
the smooth category, which settles problems of Nirenberg and
Banchoff--K\"uhnel. 
\end{abstract}

\maketitle

\section{Introduction}
A closed surface in $\R^3$ is \emph{tight} if its
total positive curvature is $4\pi$ \cite{cecil-chern1985}. These objects are natural generalizations of convex surfaces to any topological genus, and are characterized by the property that any plane cuts them into at most two pieces.
Alexandrov \cite{alexandrov1938} proved that analytic tight surfaces are \emph{rigid}:
they are congruent whenever they are isometric.
Nirenberg \cites{nirenberg1963,han-hong2006,han-hong-huang2017} obtained a smooth version under additional hypotheses, including that each component of the negatively curved region contains at most one closed
asymptotic curve.
Han and Khuri \cite{han-khuri2011} allowed several closed asymptotic
curves satisfying a nondegeneracy condition. Nirenberg had conjectured \cites{yau2000,ghomi-raffaelli2025} that
the closed asymptotic curves do not occur, which was
disproved in \cite{banks-ghomi-qiu-raffaelli2026a}.
On the other hand Banchoff \cite{banchoff:rigid} constructed nonrigid tight
polyhedral tori, and Banchoff and K\"uhnel asked whether two
isometric smooth tight tori must be congruent
\cite{banchoff-kuhnel1997}.
Here we settle these problems:

\begin{theorem}\label{thm:main}
For every integer $g\ge1$, there is a Riemannian metric on the closed
orientable surface $M$ of genus $g$ which admits uncountably many
smooth ($\C^\infty$) noncongruent tight isometric embeddings into $\R^3$.
These embeddings coincide outside a compact set in the
negatively curved region of $M$ and enclose different volumes.
Furthermore $0$ is a regular value of the curvature of $M$.
\end{theorem}

Thus Alexandrov's rigidity theorem does not hold in the smooth category, and the answer to the Banchoff--K\"uhnel question is also negative. Furthermore,  these surfaces satisfy the other hypothesis in Nirenberg's theorem that the gradient of the curvature does not vanish on the boundary of the negatively curved region; so his assumption on closed asymptotic curves cannot be removed. The first-named author and Spruck \cite{ghomi-spruck2020} showed that a positively curved surface with boundary is rigid if its boundary is held fixed; the above theorem shows that this fails for negative curvature.
We do not know, however, whether these surfaces are \emph{flexible},
i.e., admit a continuous one-parameter family of noncongruent isometric embeddings. 

The proof of Theorem~\ref{thm:main} is based on the
following observation, which goes back to
Cohn-Vossen \cite{cohn-vossen1930}; see also
Rembs \cite{rembs1952}.
Let $X\colon M\to\R^3$ be an embedding of a closed surface with induced
metric $\mathrm{I}_X\coloneqq\langle dX,dX\rangle$. A vector field $Y$ along $X$
is an \emph{infinitesimal bending} \cite{spivak:v5} if
$\langle dX(v),dY(v)\rangle=0$ for all tangent vectors $v$.
Then
\begin{equation}\label{eq:metric}
\mathrm{I}_{X+\epsilon Y}\;=\;\mathrm{I}_{X-\epsilon Y}
\;=\;\mathrm{I}_X+\epsilon^2\langle dY,dY\rangle.
\end{equation}
For small $\epsilon$, the embeddings $X\pm\epsilon Y$ therefore induce
the same metric, while their enclosed volumes differ to first order by
$2\epsilon\int\langle Y,\nu\rangle$, where $\nu$ is the outward
unit normal of $X$. Thus it suffices to
construct a tight embedding $X\colon M\to\R^3$ carrying an infinitesimal bending of nonzero
\emph{flux} $\int_M\langle Y,\nu\rangle$ supported in the negatively curved region, for then $X\pm\epsilon Y$
remain tight for small $\epsilon$.
Uncountably many embeddings are then obtained from a sequence of such
bending fields with disjoint supports.

\begin{figure}[h]
\centering
\includegraphics[width=0.75\linewidth]{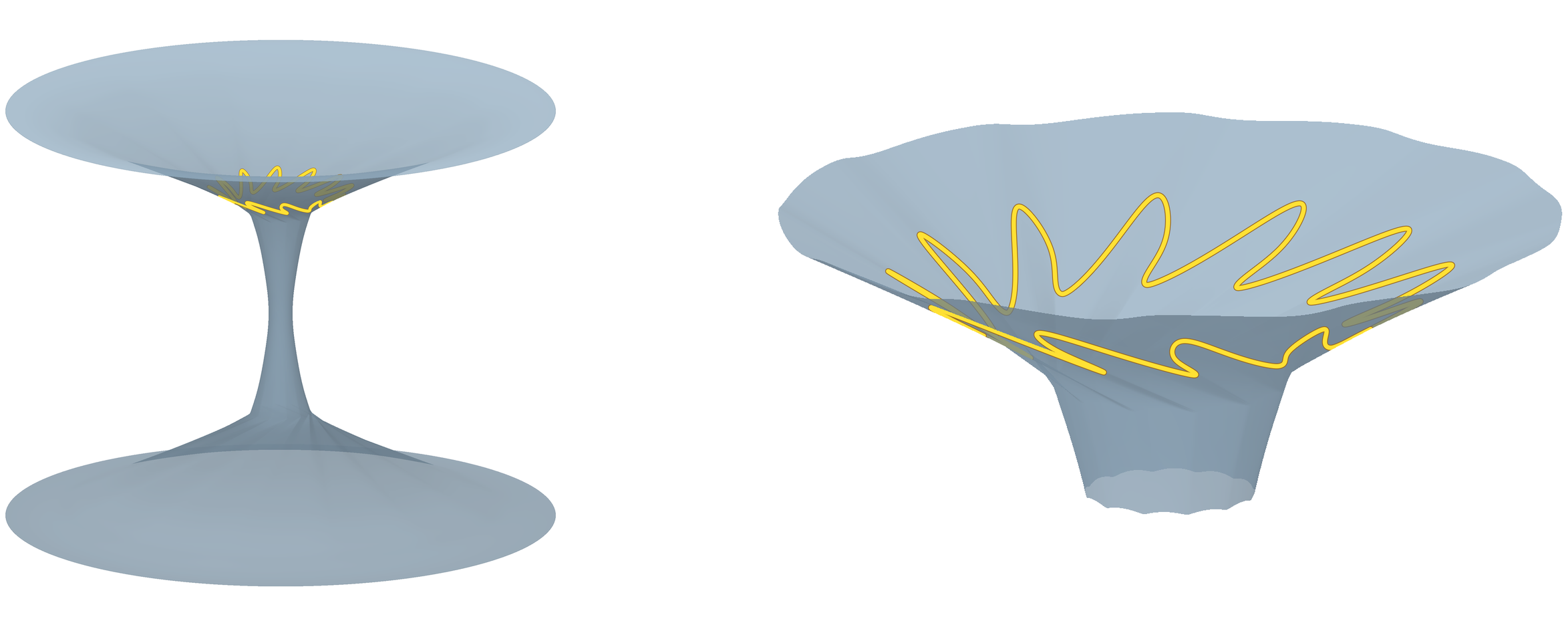}
\caption{The doubled annulus (left) and an enlarged view of its upper neck
(right). The band is shown in yellow.}
\label{fig:annulus}
\end{figure}

The above plan will be implemented as follows.
In Section~\ref{sec:band} we construct a negatively curved band ruled
by the principal normals of a closed space curve. Near this curve
the asymptotic curves transverse to the rulings are closed, which
allows us to construct infinitesimal bendings of the band supported
between any two of them. We then show that their flux can be made nonzero.
In Section~\ref{sec:completion} we extend the band to a negatively curved annulus
whose ends are surfaces of revolution.
Finally, in Section~\ref{sec:torus} we double this annulus, see Figure~\ref{fig:annulus}, 
and insert $g$ copies of it between the flat faces of a convex
surface, see Figure~\ref{fig:tight}, which yields a tight surface of genus $g$. A perturbation of this surface near the flat faces then makes $0$ a regular value of its curvature. The proofs are entirely analytic. The accompanying Mathematica notebook
\cite{ghomi-qiu-raffaelli-zhang2026} illustrates the construction with
the $12$-fold symmetric example shown in the figures in this work.

The study of rigidity of surfaces goes back to Cauchy
\cites{cauchy1813,aigner-ziegler1999} who proved it for convex polyhedra. Cohn-Vossen
\cites{cohn-vossen1927,ghomi-spruck2020} obtained the smooth version of Cauchy's
theorem, and Pogorelov \cite{pogorelov1973} finally established the
result for all convex surfaces. On the other hand, Connelly
\cite{connelly1977} constructed examples of flexible closed
polyhedral surfaces, and no smooth counterparts have yet been found.
See  \cite{bobenko-hoffmann-sageman-furnas2025} for some remarkable recent
examples of nonrigid analytic surfaces. Asymptotic curves play an
important role in rigidity of nonconvex surfaces, since they form
characteristics of the hyperbolic PDE governing the isometric embeddings
of the negatively curved regions; see
\cites{ni-zhang-zhang2026,ghomi-raffaelli2025}
for other recent results on closed asymptotic curves.

\section{Construction of the Band}\label{sec:band}

We begin by constructing a negatively
curved ruled band $X$ which, near its central curve, is foliated by
closed asymptotic curves, and
show that it carries in its interior compactly supported
infinitesimal bending fields $Y$ with nonzero flux. The band is ruled
by the principal normals of a closed space
curve $\alpha$, which is determined by its
binormal indicatrix $B$ and its speed $v$.
Sections~\ref{subsec:band}--\ref{subsec:bending} assume $B$ and
$v$ are given, and produce the band $X$ and its bending fields $Y$.
Section~\ref{subsec:flux} computes the flux of $Y$.
Section~\ref{subsec:curve} constructs $B$ and $v$, and
Section~\ref{subsec:embeddedness} shows that $X$ is embedded.

\subsection{The band $X$}\label{subsec:band}
Let $B$ be a smooth simple closed curve in the open upper hemisphere
of $\Sph^2$, parametrized by arclength $s$ with period $L$.
Set $N\coloneqq B'$, $T\coloneqq -B\times N$, and write
$N'=-B-\kappa T$, $T'=\kappa N$, which defines the geodesic curvature
$\kappa$ of $B$.
We seek a positive $L$-periodic speed function $v$ satisfying
the conditions
\begin{equation}\label{eq:conditions}
\Lambda(v)\coloneqq\int_0^L v\, T\; ds=0,
\qquad\qquad
\Pi(v)\coloneqq\int_0^L \frac{\kappa'}{\sqrt v}\, ds=0.
\end{equation}
The first condition above closes the space curve $\alpha$ defined below 
by Lelieuvre's formula \cite{lelieuvre1888}, which was used in \cite{banks-ghomi-qiu-raffaelli2026a}:
$$
\alpha'\coloneqq v T=-v B\times B'.
$$
The second condition closes the asymptotic curves on the band $X$ we describe next.
Here $(T,N,B)$ is the Darboux frame of $\alpha$ with respect to
$X$; it coincides with the Frenet frame of
$\alpha$ up to sign wherever the curvature of $\alpha$ does not
vanish, but is defined even at inflections
\cite{ghomi-raffaelli2025}*{Sec.~2}.
Define $X\colon \R/(L\mathbf{Z})\times[-\delta,\delta]\to\R^3$ by
$$
X(s,u)\coloneqq \alpha(s)+u\, N(s),
$$
where the \emph{thickness} $\delta>0$ is a small constant to be specified later in
Section~\ref{subsec:embeddedness}.
Setting $w\coloneqq v-\kappa u$ and $E\coloneqq w^2+u^2$, we compute
that $X_s=wT-uB$, $X_u=N$, and $X_s\times X_u=wB+uT$.
Hence the first fundamental form is
$E\,ds\otimes ds+du\otimes du$, the area measure is
$dA=\sqrt{E}\,ds\,du$, and
$
\nu\coloneqq (wB+uT)/\sqrt E
$
is a unit normal of $X$.
The second fundamental form with respect to $\nu$, multiplied by
$\sqrt E$, is
\begin{equation}\label{eq:II}
u(v'-\kappa' u)\,ds\otimes ds-v\,(ds\otimes du+du\otimes ds).
\end{equation}
In particular the curvature $K=-v^2/E^2<0$ as desired. Using these bands to generate closed asymptotic curves on negatively curved surfaces goes
back to Kovaleva \cite{kovaleva1968}. They feature also in the recent work of
Ni, Zhang, and Zhang \cite{ni-zhang-zhang2026}.

\subsection{The asymptotic curves}
A curve on a surface is \emph{asymptotic} if its velocity
is a null direction of the second fundamental form.
Since $X$ has negative curvature, precisely two asymptotic curves
pass through each of its points: a ruling, and a curve transverse to
the rulings. By \eqref{eq:II}, the latter curves are given by
$s\mapsto X(s,u(s))$ with $u(v'-\kappa'u)-2vu'=0$.
In particular $\alpha$ itself is asymptotic.
For $u\neq0$, the substitution $1/u$ makes this equation linear, and
multiplying by $\sqrt v$ then gives
\begin{equation}\label{eq:asymptotic}
\left(\frac{\sqrt v}{u}\right)'=\frac{\kappa'}{2\sqrt v}.
\end{equation}
By the second condition in \eqref{eq:conditions}, the right side of
\eqref{eq:asymptotic} is the derivative of a periodic function
$\omega$, which we choose so large that $\sqrt v/\omega<\delta$.
Then the function
$$
\xi(s,u)\coloneqq\frac{\sqrt{v}}{u}-\omega
$$
is constant along the asymptotic curves transverse to the rulings.
So the level sets $\xi=c$, $c>0$, are closed asymptotic curves, which
foliate the region $0<u<\sqrt v/\omega$.
By contrast, in \cite{banks-ghomi-qiu-raffaelli2026a} the rulings
were chosen first, with constant slope, and were not orthogonal to
the corresponding curve. There the analogue of
\eqref{eq:asymptotic} had only one periodic solution. On the other hand,
the surface was more explicit, and its rulings swept out an
annulus which only had to be made
tangent to planes along its boundary. 

\subsection{The bending field $Y$}\label{subsec:bending}
Now we construct the desired vector fields on the band $X$. Each is
supported between two closed asymptotic curves, and is determined
by a function $\phi$, which we call its \emph{profile function}.

\begin{proposition}\label{prop:bending}
Let $I=(c_1,c_2)\subset(0,\infty)$ be a bounded open interval, and
$\phi\in\C^\infty_c(I)$. Then the field $Y\coloneqq PT+QB$, where
\begin{equation}\label{eq:bendingfield}
P\coloneqq-\frac{u}{2\sqrt v}\,\phi'(\xi),
\qquad\qquad
Q\coloneqq u\,\phi(\xi)-\frac{w}{2\sqrt v}\,\phi'(\xi),
\end{equation}
is an infinitesimal bending of $X$, whose support lies in the region
$c_1<\xi<c_2$ in the interior of the domain of $X$.
\end{proposition}
\begin{proof}
Since $X_u=N$ and $Y_u=P_uT+Q_uB$, we have $\langle X_u,Y_u\rangle=0$.
Thus $Y$ is an infinitesimal bending if and only if
$\langle X_s,Y_s\rangle=0$ and
$\langle X_s,Y_u\rangle+\langle X_u,Y_s\rangle=0$.
Substituting $X_s=wT-uB$ and $Y_s=P_sT+Q_sB+(\kappa P+Q)N$ gives
$$
wP_s-uQ_s=0,\qquad wP_u-uQ_u+\kappa P+Q=0.
$$
Writing $Q=uf+(w/u)P$, and using $w_u=-\kappa$ and $w+\kappa u=v$,
the second equation becomes $P=u^3f_u/(2v)$, and then the first
equation becomes
$$
f_s+\frac{u(v'-\kappa'u)}{2v}\,f_u=0.
$$
Thus $f$ must be constant along the asymptotic curves transverse to the
rulings, which are the characteristics of this equation, and
conversely any such $f$ yields a bending field via the above formulas.
Since $\xi$ is constant along these curves we may take
$f\coloneqq\phi(\xi)$, and then $\xi_u=-\sqrt v/u^2$ yields
\eqref{eq:bendingfield}.
Since $\phi$ has compact support in $I$, the support of $Y$ lies in
the region $c_1<\xi<c_2$. This region lies in the interior of the
domain of $X$, since $\sqrt v/\omega<\delta$, and so $Y$ extends
smoothly by zero.
\end{proof}

\begin{figure}[h]
\centering
\begin{overpic}[width=\linewidth]{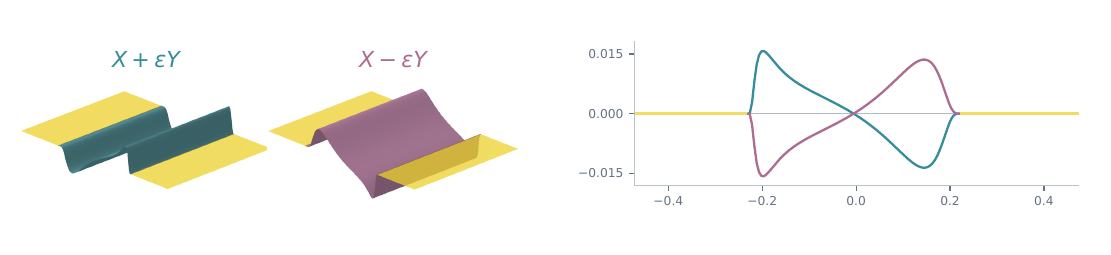}
\end{overpic}
\caption{Bendings of a patch of the ruled band, with the normal
coordinate enlarged eightfold. At right are the orthogonal projections of the curves
$u\mapsto(X\pm\epsilon Y)(s_0,u)$ onto the normal plane of $\alpha$
at $\alpha(s_0)$.}
\label{fig:bendings}
\end{figure}

Figure~\ref{fig:bendings} shows the effect of the two signs of $Y$ on
the same patch of the band.

\subsection{The flux}\label{subsec:flux}
Next we compute the flux  of the bending field $Y$
over the band, with respect to the unit normal $\nu=(wB+uT)/\sqrt E$
of $X$ from Section \ref{subsec:band}, and show that the profile
function $\phi$
can be chosen to make this flux nonzero.

\begin{lemma}\label{lem:flux}
Let $Y$ be the bending field \eqref{eq:bendingfield} with profile
$\phi\in\C_c^\infty(I)$.
Then
\begin{equation}\label{eq:kernel}
\int \langle Y,\nu\rangle
=2\int_I\phi(c)\int_0^{L}
\left(
\frac{\kappa\, v^{3/2}}{(c+\omega)^{4}}
-\frac{(1+\kappa^2)\, v}{(c+\omega)^{5}}
\right)\,ds\,dc .
\end{equation}
\end{lemma}
\begin{proof}
Since $dA=\sqrt E\,ds\,du$ and
$\sqrt E\,\nu=wB+uT$, the flux density is
$$
\big\langle Y,\,wB+uT\big\rangle
= Pu+Qw
= uw\,\phi(\xi)-\frac{E}{2\sqrt v}\,\phi'(\xi).
$$
Since $\xi_u=-\sqrt v/u^2$, we may integrate the $\phi'$ term by parts
in $u$. There is no boundary term, and the coefficient of $\phi(\xi)$ is
$$
uw-\p_u\!\left(\frac{Eu^2}{2v}\right)
=\frac{2u^2}{v}(\kappa w-u).
$$
Now put $c=\xi$. Substitution of $u=\sqrt v/(c+\omega)$ and
$|du/dc|=u^2/\sqrt v$ gives \eqref{eq:kernel}.
\end{proof}

\begin{proposition}\label{prop:nonzero}
For every nonempty open interval $I\subset(0,\infty)$ there is a
function $\phi\in\C^\infty_c(I)$ whose bending field
\eqref{eq:bendingfield} has nonzero flux.
\end{proposition}
\begin{proof}
Suppose every $\phi\in\C^\infty_c(I)$ gives zero flux.
Then the inner integral in \eqref{eq:kernel}, being continuous in $c$,
vanishes throughout $I$. At any fixed $c\in I$, its $j$th derivative
therefore vanishes for every integer $j\ge0$, which gives
$$
\int_0^{L}\frac{\kappa\, v^{3/2}}{(c+\omega)^{j+4}}\;ds
\;=\;\frac{j+4}{4}
\int_0^{L}\frac{(1+\kappa^2)\,v}{(c+\omega)^{j+5}}\;ds .
$$
Dividing by $\int_0^L(c+\omega)^{-j-4}\,ds$ gives
$(j+4)\min\big((1+\kappa^2)v/(c+\omega)\big)
\le4\max(\kappa v^{3/2})$.
The minimum is positive, so this is impossible for large $j$.
\end{proof}

\subsection{The closing conditions}\label{subsec:curve}
We first choose the spherical curve $B$ and a positive speed function $v_0$ satisfying the first
condition in \eqref{eq:conditions}. Then we perturb $v_0$ to satisfy
the second condition without losing the first. Identify $\R^3$ with $\R^2\times\R$, and call the first factor the
horizontal plane. Let $J$ denote counterclockwise rotation of $\R^2$
through $\pi/2$, and
$e_\theta\coloneqq(\cos\theta,\sin\theta)$.

\subsubsection{The first condition}\label{subsubsec:first}
\begin{figure}[h]
\centering
\includegraphics[width=0.6\linewidth]{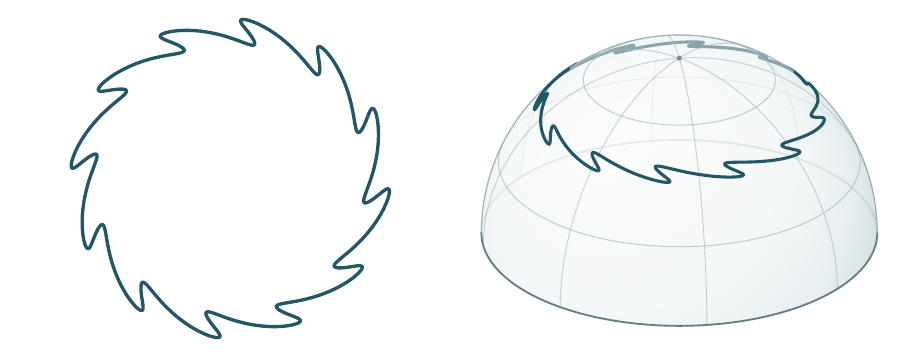}\par
\makebox[0.36\linewidth]{$\beta$}%
\makebox[0.39\linewidth]{$B$}\par
\caption{The plane curve $\beta$ with frequency $q=12$ and its spherical lift $B$.}
\label{fig:beta-B}
\end{figure}

We will obtain $B$ as the radial projection into $\Sph^2$ of a
closed curve $(\beta,1)$, where $\beta$ is a closed curve in $\R^2$;
then $B$ lies in the open upper hemisphere, as required, and is
simple whenever $\beta$ is; see Figure~\ref{fig:beta-B}.
We take $\beta$ invariant under rotation through $2\pi/q$, where the
\emph{frequency} $q\ge2$ is an integer to be chosen sufficiently large in
Section~\ref{subsec:embeddedness}, so that
the horizontal part of the first condition in \eqref{eq:conditions} holds automatically.
The vertical part then forces $\det(\beta,d\beta/d\theta)$ to change
sign, so $\beta$ cannot be star-shaped about the origin. We therefore
let the angle of $\beta$ backtrack periodically, and add a radial
wobble which separates the strands where it does, so that $\beta$
remains simple. Set
$$
\beta(\theta)\coloneqq r(\theta)e_{\widetilde\theta},\qquad
r(\theta)\coloneqq 1+\frac{\sin q\theta}{q},\qquad
\widetilde\theta\coloneqq\theta+\frac{2\sin q\theta}{q}.
$$
The curve $\beta$ is simple: if $\beta(\theta_1)=\beta(\theta_2)$, then
$r(\theta_1)=r(\theta_2)$, so $\sin q\theta_1=\sin q\theta_2$, and
the corresponding values of $\widetilde\theta$ differ by $\theta_1-\theta_2$,
which is therefore a multiple of $2\pi$.
Furthermore,
$$
\frac{d\beta}{d\theta}
=\cos q\theta\,e_{\widetilde\theta}+r(1+2\cos q\theta)\,Je_{\widetilde\theta}
$$
never vanishes. Thus $B=(\beta,1)/\sqrt{1+r^2}$ is a smooth simple
closed curve, as required in Section~\ref{subsec:band}.
Increasing $\theta$ by $2\pi/q$ rotates $\beta$ through $2\pi/q$ about
the origin, and $B$ through $2\pi/q$ about the $z$-axis.
We regard $\theta$ as a function of the arclength $s$ of $B$; since
rotations preserve arclength, increasing $s$ by $L/q$ increases
$\theta$ by $2\pi/q$.
Since $T=(dB/ds)\times B$, direct
calculation gives
$$
T\,ds=\frac{1}{1+r^2}
\left(-J\frac{d\beta}{d\theta},
-\det\left(\beta,\frac{d\beta}{d\theta}\right)\right)d\theta.
$$
For any $L/q$-periodic speed $v$, the horizontal part of
$\int_0^L vT\,ds$ vanishes by rotational symmetry. To make its vertical
part vanish as well, set
$$
\lambda(t)\coloneqq 1-\cos t+\frac14\cos 2t
=\frac14+\frac12(1-\cos t)^2,\qquad
v_0(\theta)\coloneqq\frac{1+r(\theta)^2}{r(\theta)^2}\,\lambda(q\theta).
$$
Since $\det(\beta,d\beta/d\theta)=r^2(1+2\cos q\theta)$, the vertical
component of $\int_0^L v_0T\,ds$ is
$$
-\int_0^{2\pi}\lambda(q\theta)(1+2\cos q\theta)\,d\theta=0,
$$
because the means of $\lambda(t)$ and $\lambda(t)\cos t$ are $1$ and $-1/2$,
respectively. Thus $v_0>0$ satisfies the first condition in
\eqref{eq:conditions}.

\subsubsection{The second condition}
We now decrease the speed $v_0$ on a short interval where $\kappa'$
has the sign opposite to that of $\Pi(v_0)$, which makes $\Pi$ change
sign, while a compensating perturbation elsewhere preserves the first
condition.
Set $\Lambda_3(v)\coloneqq\langle\Lambda(v),e_3\rangle$. For
$L/q$-periodic $v$ the horizontal components of $\Lambda(v)$ vanish
by symmetry, as in Section~\ref{subsubsec:first}, so the first
condition $\Lambda(v)=0$ reduces to $\Lambda_3(v)=0$.
If $\Pi(v_0)=0$, take $v=v_0$. Otherwise $\kappa$ is not constant,
and since it is $L/q$-periodic, $\kappa'$ takes both signs on a
period. Choose a short
interval on which $\kappa'$ has sign opposite to $\Pi(v_0)$, and a
smooth cutoff $0\le\chi\le1$ supported there and equal to one on a
smaller interval. Since $\langle T,e_3\rangle$ is not identically zero, we may also
choose a fixed bump $\psi$ away from the interval with
$\Lambda_3(\psi)\neq0$. Extend both cutoffs $L/q$-periodically and set
$$
v_\sigma\coloneqq v_0(1-\sigma\chi)
+\sigma\frac{\Lambda_3(v_0\chi)}{\Lambda_3(\psi)}\,\psi,
\qquad 0\le\sigma<1.
$$

The last term ensures that $\Lambda_3(v_\sigma)=0$. By shrinking the
support of $\chi$, we make its coefficient arbitrarily small. Since
the supports are disjoint, the speeds remain positive and we may
arrange that $\int_0^L|v_\sigma-v_0|\,ds<1/q$ for all $\sigma$.
Where $\chi=1$, the contribution to $\Pi(v_\sigma)$ diverges as
$\sigma\to1$ with sign opposite to $\Pi(v_0)$. The rest of the support
of $\chi$ contributes with the same sign, and the integral outside
that support stays bounded. Hence $\Pi(v_\sigma)=0$ for some $\sigma$.
Taking $v=v_\sigma$ gives both conditions in \eqref{eq:conditions}.
The curve $B$ and speed $v$ now define the band $X$ of
Section~\ref{subsec:band}.

\subsection{Embeddedness}\label{subsec:embeddedness}
Finally we show that $X$ is embedded when the frequency $q$ is
sufficiently large and the thickness $\delta$ is sufficiently small.
By compactness, it suffices to check that $\alpha$ is embedded.
 Let $\overline\alpha$ be the projection of $\alpha$ into
the horizontal plane, and translate $\alpha$ so that $\overline\alpha$
has mean zero with respect to $\theta$. We claim that, uniformly in
$\theta$,
$$
\beta(\theta)=e_\theta+O(1/q),\qquad
\overline\alpha(\theta)=\tfrac12 e_\theta+O(1/q).
$$
The first estimate follows from the definitions of $r$ and $\widetilde\theta$.
For the second, let $\overline\alpha_0$ be the projection of the curve
with speed $v_0$. By the formula for $T\,ds$,
$$
\frac{d\overline\alpha_0}{d\theta}
=-\frac{\lambda(q\theta)}{r^2}\,J\frac{d\beta}{d\theta}
=\lambda(q\theta)\big((1+2\cos q\theta)e_\theta
-\cos q\theta\,Je_\theta\big)+O(1/q).
$$
The coefficients of $e_\theta$ and $Je_\theta$ on the right have means
zero and $1/2$, respectively. A periodic function of mean zero evaluated
at $q\theta$ has an antiderivative which is $O(1/q)$. Thus, by integration
by parts, its integral against $e_\theta$ or $Je_\theta$ on any
subinterval of $[0,2\pi]$ is $O(1/q)$. Since
$d e_\theta/d\theta=Je_\theta$, it follows that
$\overline\alpha_0=\tfrac12 e_\theta+O(1/q)$ up to translation.
Furthermore, if $\overline\alpha$ and $\overline\alpha_0$ start at the
same point, then
$|\overline\alpha-\overline\alpha_0|\le\int_0^L|v-v_0|\,ds<1/q$.
The choice of mean zero now gives the claim.
In particular, $\overline\alpha$ winds once around the origin.
To see that $\overline\alpha$ is simple, note that $d\overline\alpha/d\theta$
is a positive multiple of $-J\,d\beta/d\theta$, and
$$
\Big\langle -J\frac{d\beta}{d\theta},\,e_{\theta+\pi/3}\Big\rangle
=\frac12\Big(1+\big(2-\sqrt3\,\big)\cos q\theta\Big)+O(1/q).
$$
This is bounded below by a positive constant. If
$\overline\alpha(\theta_1)=\overline\alpha(\theta_2)$ for
$0<\theta_2-\theta_1\le\pi$, the uniform estimate for $\overline\alpha$
gives $\theta_2-\theta_1=O(1/q)$. Since $d\beta/d\theta$ is uniformly
bounded, the displayed inequality remains strict on this interval
with the fixed vector $e_{\theta_1+\pi/3}$ in place of
$e_{\theta+\pi/3}$. Thus $\langle\overline\alpha,e_{\theta_1+\pi/3}\rangle$ is
strictly increasing between $\theta_1$ and $\theta_2$, a contradiction.
Thus $\overline\alpha$, and consequently $\alpha$, is embedded. So $X$
is embedded for $\delta$ sufficiently small.

\section{Extension of the Band}\label{sec:completion}
Here we extend the band $X$ to
an annulus with rotationally symmetric ends (this
surface forms the bottom half of the annulus depicted in
Figure~\ref{fig:annulus}). We work
with the support function $h$ of $X$, a function on an annular domain 
 whose gradient maps the domain onto the horizontal projection 
of $X$.
This is the classical device, going back to Pogorelov and Cheng--Yau, which
links the Minkowski problem to the Monge--Amp\`ere equation
\cites{pogorelov1978,cheng-yau1977,trudinger-wang2008}. In these
terms negative curvature means $\det\Hess h<0$, and rotational symmetry
means that $h$ is radial.

\subsection{The support function and its Legendre transform}
\label{subsec:support}
Choosing the thickness $\delta$ of $X$ sufficiently small, we may assume that
the Gauss map $\nu$ of $X$ is one-to-one into the open upper
hemisphere, and the horizontal projection of $X$ is an embedding;
indeed, along $\alpha$ these maps restrict to the embedded curves $B$
and $\overline\alpha$, and are local diffeomorphisms, since $K<0$ and
$B$ lies in the open upper hemisphere.

The \emph{normal coordinate} of a point $p$ of $X$ is the point
$x(p)\in\R^2$ which is determined by  $\nu(p)=(x(p),1)/\sqrt{1+|x(p)|^2}$. Thus $(x(p),1)$ is the
point where the ray from the origin through $\nu(p)$ meets the plane
$z=1$. In particular the normal coordinate of $\alpha$ is $\beta$, by
Section~\ref{subsubsec:first}.
Conversely, since $\nu$ is one-to-one, $x$ determines $p$, and we let
$$
h(x)\coloneqq\langle p(x),(x,1)\rangle
$$
be the \emph{support function} of $X$. Since $(x,1)$ is normal to $X$
at $p(x)$, differentiation shows that $\nabla h(x)$ is the horizontal
projection of $p(x)$. So $X$ is parametrized by
\begin{equation}\label{eq:support}
x\longmapsto\big(\nabla h,\ h-\langle x,\nabla h\rangle\big),
\end{equation}
with curvature $K=1/\big((1+|x|^2)^2\det\Hess h\big)$. In particular
$\nabla h(\beta)=\overline\alpha$, $\nabla h$ is one-to-one, and
$\det\Hess h<0$.
We say that a smooth function $h$ on a domain in $\R^2$ is
\emph{saddle} if $\det\Hess h<0$. Then, conversely,
\eqref{eq:support} is a negatively curved
immersion with normal $(x,1)$. It is an embedding if $\nabla h$ is one-to-one, in
which case we let
$$
h^*(y)\coloneqq \langle x,y\rangle-h(x)
$$
be the \emph{Legendre transform} of $h$, where $x$ is determined by
$\nabla h(x)=y$. Then $\nabla h^*(y)=x$ and
$\Hess h^*=(\Hess h)^{-1}$; in particular $\det\Hess h^*<0$.
Furthermore, by \eqref{eq:support}, the surface with support function
$h$ is the graph of $-h^*$.

\subsection{The extension problem}\label{subsec:problem}
Recall from Proposition~\ref{prop:bending} that the bending field $Y$
is supported in the region $c_1<\xi<c_2$ of $X$. By
Proposition~\ref{prop:nonzero}, we may assume that $Y$ has nonzero
flux.
Let $\Omega$ be the domain of $h$, and $\Omega^*$ be the horizontal
projection of $X$. Then 
$
\nabla h\colon \Omega\to\Omega^*
$ 
is a diffeomorphism, and $X$ is the graph of $-h^*$ over $\Omega^*$.
Denote the normal coordinates of the boundary curves of $X$, in
order from the inside, by $\gamma_1(s)$ and $\gamma_2(s)$. Then $\Omega$
is the annulus bounded by $\gamma_1$ and $\gamma_2$, and $\Omega^*$ is the annulus
bounded by the curves $\gamma_i^*\coloneqq\nabla h(\gamma_i)$.
These four curves wind counterclockwise about the origin, while
$\nabla h$ reverses orientation, since $\det\Hess h<0$. So
$\gamma_2^*$ is the inner boundary of $\Omega^*$.
To establish the desired extension of $X$ it suffices to show:

\begin{proposition}\label{prop:potential}
After a vertical translation of $X$, and a modification in an arbitrarily small collar of  $\partial\Omega$, the support function $h$ of $X$ extends to
a saddle function on $\R^2\setminus\{0\}$, so that $h$
is radial near the origin, and $h(x)=a|x|-b/|x|$ near
infinity, for $a$, $b>0$.
Moreover $\nabla h$ is a diffeomorphism of $\R^2\setminus\{0\}$ onto
an annulus $a<|y|<R_0$, and $h^*$ extends smoothly to  a constant function
across the circle $|y|=R_0$.
\end{proposition}

The extended $h$ is the support function of
an embedded negatively curved annulus $A$, which agrees with $X$ on
the bending region $c_1\le\xi\le c_2$, and is the graph of $-h^*$ over the
annulus $a<|y|<R_0$. The outer end of $A$ agrees to infinite order
with a horizontal plane, and the inner end is the surface of
revolution $|y|=a+z^2/(4b)$.
The proof is given in Section~\ref{subsec:assembly}. It rests on two
lemmas for extending a saddle function inward: across a circle, by a
radial function (Section~\ref{subsec:filling}), and across a closed
curve up to a circle, by a function which is affine along line
segments (Section~\ref{subsec:connector}). Applied to $h$, these
lemmas extend it across the inner boundary $\gamma_1$ of $\Omega$;
applied to $h^*$, they extend $h$ across the outer boundary $\gamma_2$.
We need the following fact.

\begin{lemma}[\cite{ghomi2002}*{Lem.~1.1}]
\label{lem:degree}
Let $f$ be a local diffeomorphism from a
compact connected surface with boundary into $\R^2$,
which is one-to-one on each boundary
component. Then $f$ is an embedding.
\end{lemma}

\subsection{Radial extension}\label{subsec:filling}
Here we show that a saddle function extends
radially inward across a circle, provided that its Hessian is positive
in the direction tangent to the circle (Lemma~\ref{lem:radial}); the
corresponding end of the surface is then a surface of revolution.
First we show how such functions may be joined.

\begin{lemma}\label{lem:join}
Let $h_1$ and $h_2$ be saddle functions on the two sides of a
smooth closed curve $\Gamma$ in
$\R^2$, which agree to first order along $\Gamma$.
Then within any prescribed collar $U$ of $\Gamma$ they can be joined
into a single saddle function, which
agrees with $h_1$ and $h_2$ outside $U$.
\end{lemma}
\begin{proof}
Let $f$ be the function equal to $h_1$ and $h_2$ on the two sides of
$\Gamma$, which is $\C^1$ since $h_1$ and $h_2$ agree to first order. In a
tangent--normal frame along $\Gamma$ the first columns of the
one-sided Hessians agree, by differentiation of the common gradient
along $\Gamma$. Thus the determinant is affine on the segment between
the two Hessians and is strictly negative throughout, uniformly along
the compact curve $\Gamma$. Let $f_\epsilon$ be the mollification of $f$
with radius $\epsilon$. Since $f$ is $\C^1$, its distributional Hessian
has no singular part along $\Gamma$, so near $\Gamma$ the Hessian of
$f_\epsilon$ is uniformly close to that segment, while away from $\Gamma$
we have $f_\epsilon\to f$ in $\C^2$. Fix a cutoff function $\chi$ supported in $U$
and equal to one near $\Gamma$, and put
$\widetilde f\coloneqq f+\chi(f_\epsilon-f)$. Then $\widetilde f=f_\epsilon$ near
$\Gamma$, $\widetilde f=f$ outside $U$, and on the support of
$\nabla\chi$, which lies at a positive distance from $\Gamma$,
$\widetilde f\to f$ in $\C^2$. Hence $\det\Hess\widetilde f<0$ for small
$\epsilon$.
\end{proof}

We write $\Hess f>0$ \emph{along the circle} $|y|=R$ if
$\Hess f(v,v)>0$ for all nonzero vectors $v$ tangent to
the circle, or equivalently, 
$f_{\theta\theta}+Rf_\rho>0$ on the circle in polar coordinates.

\begin{lemma}\label{lem:radial}
Let $f$ be a saddle function on a
neighborhood of the circle $|y|=R$ in $\{|y|\ge R\}$, with
$\Hess f>0$ along the circle,
and $g$ be a smooth function on $(a,R]$, $0\le a<R$, with $g'(R)=0$
and $g''<0$.
Then, after a modification in an arbitrarily small collar of the
circle, $f$ extends to a saddle function on $a<|y|\le R$, which
coincides with $\mu g(|y|)$ near $|y|=a$, up to an additive constant,
for a large constant $\mu$.
\end{lemma}
\begin{proof}
We may assume that $g(R)=0$. Choose a cutoff $\chi$ on $(a,R]$ with
$\chi=0$ near $a$ and $\chi=1$ near $R$, and put
$$
f_\mu\coloneqq\mu\,g(\rho)
+\chi(\rho)\big(f(Re_\theta)+(\rho-R)\,f_\rho(Re_\theta)\big),
$$
which agrees with $f$ to first order along the circle.
The radial and tangential entries of the Hessian of $f_\mu$ are
$\p_\rho^2f_\mu$ and
$(\p_\theta^2f_\mu+\rho\,\p_\rho f_\mu)/\rho^2$.
For large $\mu$ the first is negative and the second is positive:
the terms multiplying $\mu$ are $g''<0$ and $g'/\rho$, which is
positive for $\rho<R$, while the remaining terms are bounded, vanish
near $a$, and in the second entry equal the positive function
$(f_{\theta\theta}+Rf_\rho)/R^2$ at $\rho=R$.
So $\det\Hess f_\mu<0$ regardless of the mixed entry, and
Lemma~\ref{lem:join} joins $f$ and $f_\mu$ in a collar of the circle.
\end{proof}

\subsection{Extension to a circle}
\label{subsec:connector}
Here we extend a saddle function $f$
inward, across a closed curve $\gamma$ enclosing the origin up to a circle centered at the origin, so that
$\Hess f>0$ along that circle
(Lemma~\ref{lem:connector}). The extension is affine along line
segments which join $\gamma$ to the circle, so its graph is a ruled
surface. We need $\gamma$ and its
image $\gamma^*$ under $\nabla f$ to satisfy the following condition.

\begin{figure}[h]
\centering
\begin{overpic}[width=0.4\linewidth]{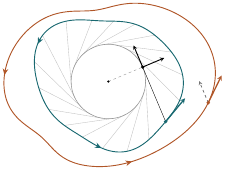}
\put(54,46){\makebox(0,0){$R$}}
\put(75,46){\makebox(0,0){$e_\varphi$}}
\put(59.0,58.5){\makebox(0,0){$Je_\varphi$}}
\put(75.4,18.1){\makebox(0,0){$\gamma$}}
\put(84.5,32.5){\makebox(0,0){$\dot\gamma$}}
\put(97,28){\makebox(0,0){$\gamma^*$}}
\put(101.5,43){\makebox(0,0){$\dot\gamma^*$}}
\put(88.5,42.5){\makebox(0,0){$Je_\varphi$}}
\end{overpic}
\caption{A pair of curves $(\gamma,\gamma^*)$ satisfying the circle
condition at $R$: $\dot\gamma$ has a positive component along
$e_\varphi$, and $\dot\gamma^*$ along $Je_\varphi$.}
\label{fig:circle}
\end{figure}
For every point $y\in\R^2$ with $|y|>R>0$ there is a unique angle
$\varphi(y)$, mod $2\pi$, such that
$
y=Re_\varphi-\ell\,Je_\varphi,
$
$\ell>0$.
So $y$ lies on the line tangent to the circle of radius $R$ about the
origin at $Re_\varphi$, and the segment from $y$ to $Re_\varphi$ points
in the counterclockwise direction $Je_\varphi$; see Figure~\ref{fig:circle}.
We say that a pair $(\gamma,\gamma^*)$ of regular closed curves in $\R^2$,
with a common parameter $s$, satisfies the \emph{circle condition} at $R$
if $|\gamma|>R$ and, 
$$
\left\langle\dot\gamma(s),e_{\varphi(\gamma(s))}\right\rangle>0,\qquad
\left\langle\dot\gamma^*(s),Je_{\varphi(\gamma(s))}\right\rangle>0.
$$
The first inequality means that $\gamma$ crosses the tangent segments
to the circle transversely, with the point of tangency $Re_\varphi$
moving counterclockwise. The second means that $\dot\gamma^*(s)$ makes an
acute angle with the segment from $\gamma(s)$ to its point of
tangency $Re_\varphi$.

\begin{lemma}\label{lem:connector}
Let $\gamma$ be a positively oriented simple closed curve in $\R^2$
enclosing the origin, $f$ be a saddle function
on a neighborhood of $\gamma$ in the closure of its
exterior, and $\gamma^*\coloneqq\nabla f(\gamma)$.
Suppose that $(\gamma,\gamma^*)$ satisfies the circle condition at $R$.
Then, after a modification in an arbitrarily small collar of
$\gamma$, $f$ extends to a saddle function on the region
between $\gamma$ and the circle $|y|=R$, so that $\Hess f>0$ along
this circle.
\end{lemma}
\begin{proof}
Put $\varphi\coloneqq\varphi(\gamma)$. Then
$\gamma=Re_\varphi-\ell\,Je_\varphi$, where
$\ell\coloneqq\sqrt{|\gamma|^2-R^2}$, and differentiation yields
$\dot\varphi=\langle\dot\gamma,e_\varphi\rangle/\ell>0$.
Since $\gamma$ winds once about the origin, $\varphi$ increases by $2\pi$.
For a small constant $\tau>0$ set $\theta\coloneqq\varphi-\tau$ and
$W\coloneqq Re_\theta-\gamma$.
The quantities
$$
D_1\coloneqq\langle\dot\gamma^*,W\rangle,
\qquad
D_2\coloneqq\det(\dot\gamma,W),
\qquad
D_3\coloneqq\det\big(R\,\dot\theta\,Je_\theta,W\big)
=R\,\dot\theta\,\big(\langle\gamma,e_\theta\rangle-R\big)
$$
are positive for sufficiently small $\tau$: since
$W\to\ell Je_\varphi$ as $\tau\to0$, the first two tend
to $\ell\langle\dot\gamma^*,Je_\varphi\rangle$ and $\ell\langle\dot\gamma,e_\varphi\rangle$,
while
$\langle\gamma,e_\theta\rangle-R
=\ell\sin\tau-R(1-\cos\tau)>0$.
The map $(s,t)\mapsto\gamma+tW$, $0\le t\le1$, has Jacobian
$\Delta\coloneqq(1-t)D_2+tD_3>0$, and is one-to-one on the boundary
curves $\gamma$ and $Re_\theta$. So, by Lemma~\ref{lem:degree}, it is
a diffeomorphism onto the region between $\gamma$ and the circle
$|y|=R$. On this region let $\widehat f$ be the function which is affine
on each segment $t\mapsto\gamma+tW$, and agrees with $f$ to first
order at $\gamma$, i.e.,
$$
\widehat f(\gamma+tW)\coloneqq f(\gamma)+t\langle\gamma^*,W\rangle.
$$
Differentiating with respect to $s$ and $t$ yields
$
\nabla\widehat f(\gamma+tW)=\gamma^*-\frac{tD_1}{\Delta}\,JW.
$
In particular $\widehat f$ agrees with $f$ to first order along $\gamma$.
Since $\widehat f$ is affine along the segments, $W$ is a null direction
of $\Hess\widehat f$. So $\det\Hess\widehat f\le0$, with equality only where
$(\Hess\widehat f)W=\p_t\nabla\widehat f$ vanishes; but
$\p_t\nabla\widehat f=-(D_1D_2/\Delta^2)\,JW\neq0$.

It remains to check that $\Hess\widehat f>0$ along the circle. There
$\widehat f(Re_\theta)=f(\gamma)+\langle\gamma^*,W\rangle$, so
$\widehat f_{\theta\theta}$ is bounded uniformly as $\tau\to0$, since
$\dot\varphi>0$. On the other hand, by the above formula for
$\nabla\widehat f$ with $t=1$,
$\widehat f_\rho=\langle\gamma^*,e_\theta\rangle
+(D_1/D_3)\langle W,Je_\theta\rangle$.
As $\tau\to0$, $D_3\to0$ while $D_1$ and $\langle W,Je_\theta\rangle$
tend to the positive quantities
$\ell\langle\dot\gamma^*,Je_\varphi\rangle$ and $\ell$. So
$\widehat f_\rho\to\infty$ uniformly, and
$\widehat f_{\theta\theta}+R\widehat f_\rho>0$ for small $\tau$.
Finally, Lemma~\ref{lem:join} joins $f$ and $\widehat f$ in a collar of
$\gamma$.
\end{proof}

\subsection{Proof of Proposition~\ref{prop:potential}}\label{subsec:assembly}
First we check that the boundary curves of $\Omega^*$ and $\Omega$
satisfy the circle condition, the latter after a reflection.
Let $\mathcal R$ denote reflection in
the first coordinate axis. For a curve $\gamma(s)$ we also write
$\mathcal R\gamma$ for its reflection with the parameter reversed,
i.e., the curve $s\mapsto\mathcal R(\gamma(-s))$; so $\mathcal R\gamma$
winds counterclockwise whenever $\gamma$ does.

\begin{lemma}\label{lem:circle}
If the frequency $q$ of $\beta$ is sufficiently large, and the
thickness $\delta$ of $X$ is sufficiently small, then
$(\gamma_2^*,\gamma_2)$ and $(\mathcal R\gamma_1,\mathcal R\gamma_1^*)$
satisfy the circle condition.
\end{lemma}
\begin{proof}
 As $\delta\to0$, $\gamma_i$ converge to $\beta$
in $\C^1$, and so $\gamma_i^*=\nabla h(\gamma_i)$ converge to
$\nabla h(\beta)=\overline\alpha$. So it
suffices to check the conditions for the pairs $(\overline\alpha,\beta)$
and $(\mathcal R\beta,\mathcal R\overline\alpha)$,
for which, by Section~\ref{subsec:embeddedness},
$\beta=e_\theta+O(1/q)$ and $\overline\alpha=\tfrac12 e_\theta+O(1/q)$
uniformly, and $d\overline\alpha/d\theta$ is a positive multiple of
$-J\,d\beta/d\theta$. In both pairs $(\gamma,\gamma^*)$, $\dot\gamma^*$ is therefore a positive
multiple of $J\dot\gamma$, so the two inequalities of the circle
condition reduce to $\langle\dot\gamma,e_\varphi\rangle>0$.
For the first pair, with $R=1/4$,
$e_\varphi=e_{\theta+\pi/3}+O(1/q)$, and
$\langle\dot\gamma,e_\varphi\rangle$ is a positive multiple of
$\langle -J\,d\beta/d\theta,\,e_{\theta+\pi/3}\rangle$, which is bounded below by
a positive constant by Section~\ref{subsec:embeddedness}.
For the second pair, $\gamma=\mathcal R\beta=e_\theta+O(1/q)$,
so, with $R=4/5$, $e_\varphi=\tfrac45 e_\theta+\tfrac35 Je_\theta+O(1/q)$,
while
$$
\dot\gamma=-\mathcal R\,(d\beta/d\theta)(-\theta)
=-\cos q\theta\,e_\theta+(1+2\cos q\theta)\,Je_\theta+O(1/q),
$$
and hence $\langle\dot\gamma,e_\varphi\rangle=(3+2\cos q\theta)/5+O(1/q)>0$.
The radius inequalities follow from the same approximations.
\end{proof}

Now we prove Proposition~\ref{prop:potential}, by extending $h$ first
inward, across the inner boundary $\gamma_1$ of $\Omega$, and then
outward, across the outer boundary $\gamma_2$.

\subsubsection{Inward extension}\label{subsubsec:inward}
By Lemma~\ref{lem:circle}, the pair
$(\mathcal R\gamma_1,\mathcal R\gamma_1^*)$ of reflected curves
satisfies the circle condition at some $R$.
Here $\mathcal R\gamma_1$ is the inner boundary of the domain
$\mathcal R\Omega$ of the saddle function $h(\mathcal Rx)$, and is
mapped by the gradient of this function onto $\mathcal R\gamma_1^*$.
So Lemma~\ref{lem:connector} extends $h(\mathcal Rx)$, after a
modification in a small collar of $\mathcal R\gamma_1$, to a saddle
function on the region between $\mathcal R\gamma_1$ and the circle
$|x|=R$, whose Hessian is positive along this
circle. This is the hypothesis of Lemma~\ref{lem:radial}, which, with
$a=0$, extends the function further to the region $0<|x|\le R$, once
we choose a function $g$ on $(0,R]$ with $g'(R)=0$ and $g''<0$.
Let $\psi\colon[0,1]\to[0,R]$ be a smooth function with $\psi(0)=R$,
$\psi(1)=0$, and $\psi'<0$ on $[0,1)$, which is flat at $1$, and
$g$ be the antiderivative of $\psi^{-1}$ on $(0,R]$ with $g(0)=0$.
Reflecting back, we then obtain an extension of $h$ to a saddle
function on the region $|x|>0$ inside $\gamma_2$, so that
$h(x)-\mu g(|x|)$ is constant near the origin, for a large constant
$\mu$. Set $R_0\coloneqq\mu$.

\subsubsection{Outward extension}
As $\nabla h$ maps the outer boundary
$\gamma_2$ of $\Omega$ onto the inner boundary $\gamma_2^*$ of $\Omega^*$, we
extend $h^*$ inward across $\gamma_2^*$, and then take the Legendre
transform. The function $h^*$ is saddle on $\Omega^*$, and its
gradient maps $\gamma_2^*$ onto $\gamma_2$. By
Lemma~\ref{lem:circle}, the pair $(\gamma_2^*,\gamma_2)$ satisfies
the circle condition at some $R$.
So Lemma~\ref{lem:connector} extends $h^*$, after a modification in a
small collar of $\gamma_2^*$, to a saddle function on the region
between $\gamma_2^*$ and the circle $|y|=R$, whose Hessian is positive
along this circle. 

This is again the hypothesis
of Lemma~\ref{lem:radial}, which we apply with $0<a<\min\{R,R_0\}$,
and a function $g$ on $(a,R]$ which equals $2\sqrt{\rho-a}$ near $a$,
and satisfies $g'(R)=0$ and $g''<0$. This extends $h^*$ to a saddle function
on the region $|y|>a$ inside $\gamma_1^*$, so that
$h^*(y)=2\sqrt{b(|y|-a)}-c$ near $|y|=a$, for some constants
$b>0$ and $c$. Since a vertical translation of $X$ changes $h$
only by an additive constant, we may assume that $c=0$.
Since the above modifications take place in small collars of
$\gamma_1$ and $\gamma_2^*$, they do not affect $h$ near $\beta$, or
near the region $c_1\le\xi\le c_2$.

On the annulus between $\overline\alpha$ and a circle $|y|=\rho$, with $\rho$
close to $a$, the gradient of $h^*$ is a local diffeomorphism which
maps the boundary curves one-to-one onto $\beta$ and the circle
$|x|=\sqrt{b/(\rho-a)}$. So it is an embedding by
Lemma~\ref{lem:degree}. Letting $\rho\to a$, we see that $\nabla h^*$
maps the region $|y|>a$ inside $\overline\alpha$ diffeomorphically onto the
exterior of $\beta$. Hence the Legendre transform of $h^*$ on this
region is a saddle function on the
exterior of $\beta$, which coincides with $h$ near $\beta$, and
equals $a|x|-b/|x|$ near infinity. Together with the
inward extension of $h$, this yields the desired extension of
$h$ to $\R^2\setminus\{0\}$.

Finally, on sufficiently small and large circles $|x|=\rho$,
$\nabla h$ equals $\mu g'(\rho)e_\theta$ and $(a+b/\rho^2)e_\theta$
respectively, and $\mu g'(\rho)\to R_0$ as $\rho\to0$. So
Lemma~\ref{lem:degree}, applied to the annuli between these circles,
shows that $\nabla h$ is a diffeomorphism onto $a<|y|<R_0$.
Furthermore, for $|y|$ close to $R_0$, the radial derivative of $h^*$
is $|x|=\psi(|y|/R_0)$, and $h^*$ converges to a constant as $|y|\to R_0$.
Since $\psi$ is flat at $1$, it follows that $h^*$ extends smoothly
by this constant across the circle $|y|=R_0$.\qed

\section{Proof of Theorem~\ref{thm:main}}\label{sec:torus}

We first assemble a tight surface $\Sigma$ of genus $g$ from copies
of the annulus $A$, and perturb it so that $0$ becomes a regular
value of its curvature. Then we bend $\Sigma$ along the fields $Y$
to obtain the desired embeddings.

\subsection{The initial surface}
\begin{figure}[h]
\centering
\includegraphics[width=0.75\linewidth]{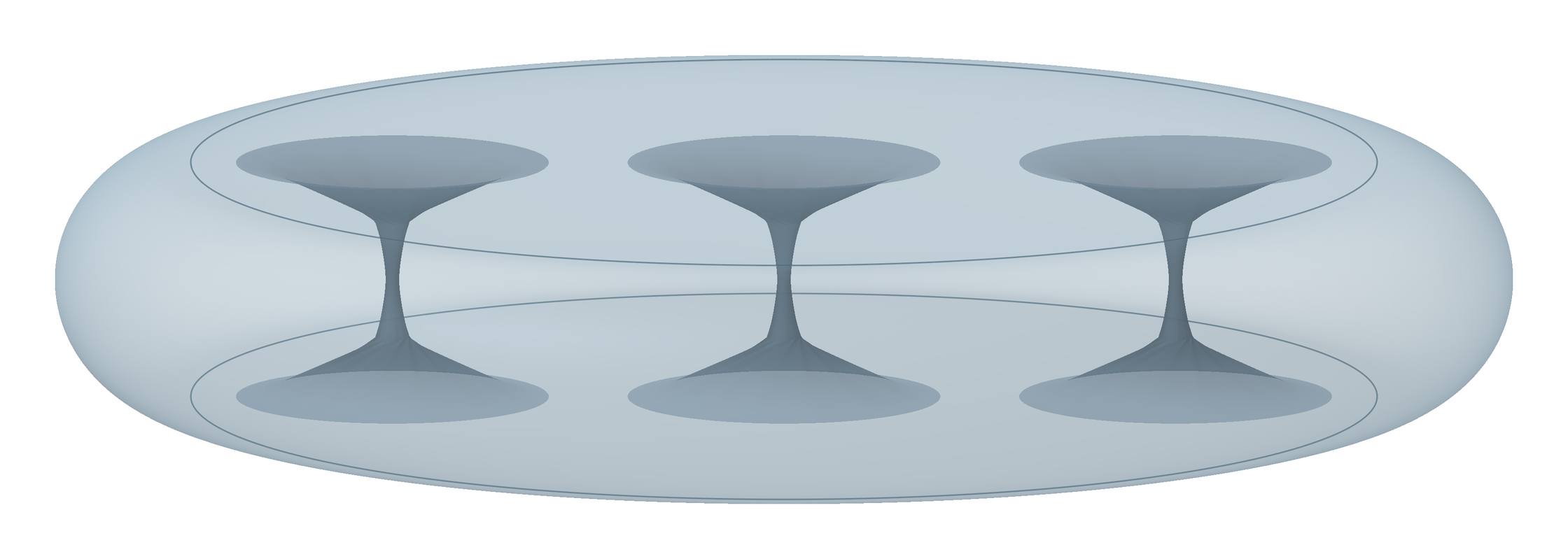}
\caption{A genus $3$ tight surface constructed using the annulus in Figure~\ref{fig:annulus}.}
\label{fig:tight}
\end{figure}

Let $A$ be the annulus constructed in Section~\ref{subsec:problem}. Recall that $A$ is the graph of $-h^*$ which is saddle, and so has no local extrema. Thus $A$ lies 
between the horizontal planes of its boundaries.
Since the top end of $A$ is the surface of revolution
$|y|=a+z^2/(4b)$, it extends smoothly by reflection in the plane $z=0$. Thus $A$ and
its reflection form a smooth embedded annulus $\widetilde A$ with
negative curvature in its interior, which agrees to infinite order with the planes of its boundary circles. Let $C$ be a smooth closed convex surface, symmetric with respect to a horizontal plane, with flat regions on top and bottom as depicted in Figure~\ref{fig:tight}. We may assume that the flat regions are large enough to contain $g$ disjoint disks of radius equal to that of the boundaries of $\widetilde A$, and that the flat regions are separated by the same distance as the height of $\widetilde A$. Removing the disks and gluing in $g$ copies of $\widetilde A$ yields a smooth closed surface $\Sigma$ of genus $g$ which is tight, since by construction its total positive curvature agrees with that of $C$ which is $4\pi$.

\subsection{The perturbation}
Next we perturb $\Sigma$ near its flat regions, so that $0$ becomes a
regular value of its curvature.
Let $z=H$ be the plane of the top flat region of $\Sigma$, and
$(y_i,H)$ be the centers of the upper boundary circles of the copies
of $\widetilde A$ in $\Sigma$. Near the $i$-th circle, $\Sigma$ is the
graph $z=H+f(r)$, where $r\coloneqq|y-y_i|$. By the last paragraph of
Section~\ref{subsec:assembly}, $f'>0$ and $f''<0$ for $r<R_0$, while
$f=0$ for $r\ge R_0$.

Now, for a small constant $\lambda>0$, choose $C$ to be positively
curved, with its top region given by the graph $z=H-\lambda|y|^2/2$
instead of being flat, and subtract
$\lambda\eta(r)|y|^2/2$ from the height of the $i$-th copy of
$\widetilde A$, where $\eta$ is a smooth function on $(0,R_0]$ which
is equal to $1$ near $R_0$, and vanishes outside a slightly larger
neighborhood of $R_0$. Modify the bottom of $\Sigma$ symmetrically.
Where $\eta$ is not constant, the curvature of $\widetilde A$ is
negative and bounded away from zero; so it remains negative for small
$\lambda$. Furthermore, the modified part of $\widetilde A$ is disjoint
from the band $X$.

Near the $i$-th circle, the new surface is the graph
$z=H+f(r)-\lambda|y|^2/2$, whose Hessian has eigenvalues $f''-\lambda<0$ and
$f'/r-\lambda$. Since $f'/r$ decreases strictly to zero as $r\to R_0$,
there is a unique radius $r_\lambda$ with $f'(r_\lambda)=\lambda r_\lambda$, which
is a simple zero of $f'/r-\lambda$. So the curvature is negative for
$r<r_\lambda$, positive for $r>r_\lambda$, and its gradient does not vanish
on the curve $r=r_\lambda$. Along this curve the gradient of the height function is $-\lambda y_i$.
So the curve lies in the plane
$$
z=H+f(r_\lambda)-\frac{\lambda}{2}r_\lambda^2+\frac{\lambda}{2}|y_i|^2
-\lambda\langle y_i,y\rangle,
$$
which is tangent to the surface along the curve. These planes have
normals $(\lambda y_i,1)$, and so are pairwise nonparallel. Since
the horizontal projection of the curve is a circle, the curve is an
ellipse, whose geodesic curvature, with respect to the
conormal which points into the negatively curved region, coincides
with its curvature in this plane. So each component of the negatively
curved region is an annulus whose boundary curves each have total
geodesic curvature $2\pi$, and therefore has total curvature $-4\pi$
by the Gauss--Bonnet theorem. Since $C$ is positively curved, the
curvature is positive elsewhere; so $0$ is a regular value of the
curvature. It follows that the total positive curvature
of the new surface is $4\pi(1-g)+4\pi g=4\pi$. So the new
surface, which we again denote by $\Sigma$, is tight.

\subsection{The bendings}
Let $X\colon M\to\R^3$ be the embedding with image $\Sigma$
constructed above. On one of the copies of $\widetilde A$ take the bending field $Y$ of
nonzero flux chosen in Section~\ref{subsec:problem}. Its support is compact in the region where $A$ agrees with the band
of Section~\ref{sec:band}, so it extends by zero to a smooth
infinitesimal bending of $X$, with flux
$\mathcal F\coloneqq\int_M\langle Y,\nu\rangle\neq0$.
Since $M$ is compact, its embeddings are open in the $\C^1$-topology.
Thus $X_\pm\coloneqq X\pm\epsilon Y$ are embeddings for small
$\epsilon>0$, and by \eqref{eq:metric} they induce the same metric on $M$.
Since $\operatorname{supp}Y$ is compact in the region of negative
curvature, $\C^2$ closeness keeps the curvature of $X_\pm$ negative
there. Hence the total positive curvature remains at $4\pi$ and tightness is preserved.
By the first variation of volume, the volumes $V(X_\pm)$ enclosed by
$X_\pm$ satisfy $V(X_+)-V(X_-)=2\epsilon\mathcal F+o(\epsilon)$.
Since $\mathcal F\neq0$, these volumes differ for small $\epsilon>0$,
and hence $X_\pm$ are noncongruent.

To obtain uncountably many embeddings, apply
Proposition~\ref{prop:nonzero} to a sequence of disjoint subintervals
of $I$. This gives bending fields $Y_1,Y_2,\dots$ of $X$ with disjoint
supports and nonzero flux. Then
$2\delta_j\coloneqq V(X+t_jY_j)-V(X-t_jY_j)$ is nonzero for small
$t_j>0$, as shown above. Choose $t_j$ successively so small that
$\delta_j\neq0$, $|\delta_j|<|\delta_{j-1}|/4$ for $j\ge2$, and
$|t_jY_j|_{\C^{j+2}}\le\epsilon2^{-j}$, where $\epsilon>0$ is a small
constant. For $\sigma\in\{-1,1\}^{\mathbf N}$ set
$$
X_\sigma\coloneqq X+\sum\sigma_jt_jY_j.
$$
 This series converges in
every $\C^k$ norm, and $|X_\sigma-X|_{\C^2}\le\epsilon$. So $X_\sigma$
is a smooth embedding, which is tight as above, and $0$ is a regular
value of its curvature, since $X_\sigma=X$ outside a compact set where
the curvature is negative. Furthermore, since at
each point of $M$ at most one $dY_j$ is nonzero, \eqref{eq:metric}
shows that $X_\sigma$ induces the metric
$\mathrm I_X+\sum_jt_j^2\langle dY_j,dY_j\rangle$, which does not
depend on $\sigma$. It remains to check that the embeddings $X_\sigma$ enclose different
volumes. Since the supports of $Y_j$ are disjoint, and
$V(X)=\frac13\int_M\langle X,\nu\rangle$ by the divergence theorem,
with the same formula for each $X_\sigma$, the contributions of $Y_j$
to the enclosed volume add. So, first for
finite sums and then in the limit by $\C^1$ convergence,
$V(X_\sigma)-V(X_{\sigma'})=\sum_j(\sigma_j-\sigma_j')\delta_j$.
If $\sigma\neq\sigma'$, and $j$ is the first index with
$\sigma_j\neq\sigma_j'$, then the absolute value of this sum is at
least $2|\delta_j|-2\sum_{k>j}|\delta_k|>0$.
\qed

\section*{Acknowledgement}
Initial drafts of the main proofs were developed with the assistance of Claude (Anthropic) and ChatGPT (OpenAI), building on our earlier work \cite{banks-ghomi-qiu-raffaelli2026a}. These drafts required extensive clarification, revision, and verification before reaching their present form. We take full responsibility for the final manuscript.

\bibliography{references}

\begin{bibdiv}
\begin{biblist}

\bib{aigner-ziegler1999}{book}{
      author={Aigner, Martin},
      author={Ziegler, G\"{u}nter~M.},
       title={Proofs from {T}he {B}ook},
   publisher={Springer-Verlag, Berlin},
        date={1999},
        ISBN={3-540-63698-6},
        note={Including illustrations by Karl H. Hofmann, Corrected reprint of
  the 1998 original},
      review={\MR{1723092}},
}

\bib{alexandrov1938}{article}{
      author={Alexandrov, A.~D.},
       title={On a class of closed surfaces},
        date={1938},
     journal={Mat. Sbornik},
      volume={4},
       pages={69\ndash 77},
}

\bib{banchoff:rigid}{article}{
      author={Banchoff, Thomas~F.},
       title={Non-rigidity theorems for tight polyhedra},
        date={1970},
     journal={Arch. Math. (Basel)},
      volume={21},
       pages={416\ndash 423},
      review={\MR{0273548}},
}

\bib{banchoff-kuhnel1997}{incollection}{
      author={Banchoff, Thomas~F.},
      author={K{\"u}hnel, Wolfgang},
       title={Tight submanifolds, smooth and polyhedral},
        date={1997},
   booktitle={Tight and taut submanifolds ({B}erkeley, {CA}, 1994)},
      series={Math. Sci. Res. Inst. Publ.},
      volume={32},
   publisher={Cambridge Univ. Press},
     address={Cambridge},
       pages={51\ndash 118},
}

\bib{banks-ghomi-qiu-raffaelli2026a}{unpublished}{
      author={Banks, Harrison},
      author={Ghomi, Mohammad},
      author={Qiu, Hongda},
      author={Raffaelli, Matteo},
       title={A counterexample to {N}irenberg's conjecture on closed asymptotic
  curves},
        date={2026},
        note={arXiv:2610.13848},
}

\bib{bobenko-hoffmann-sageman-furnas2025}{article}{
      author={Bobenko, Alexander~I.},
      author={Hoffmann, Tim},
      author={Sageman-Furnas, Andrew~O.},
       title={Compact {B}onnet pairs: isometric tori with the same curvatures},
        date={2025},
     journal={Publ. Math. Inst. Hautes \'Etudes Sci.},
      volume={142},
       pages={241\ndash 293},
}

\bib{cauchy1813}{article}{
      author={Cauchy, A.~L.},
       title={Sur les polygones et les poly\`edres, {S}econd {M}\'emoire},
        date={1813},
     journal={J. \'Ecole Polytechnique},
      volume={9},
       pages={87\ndash 98},
}

\bib{cecil-chern1985}{book}{
      editor={Cecil, Thomas~E.},
      editor={Chern, Shiing-Shen},
       title={Tight and taut submanifolds},
      series={Mathematical Sciences Research Institute Publications},
   publisher={Cambridge University Press},
     address={Cambridge},
        date={1997},
      volume={32},
        ISBN={0-521-62047-3},
      review={\MR{1486867}},
}

\bib{cheng-yau1977}{article}{
      author={Cheng, Shiu~Yuen},
      author={Yau, Shing~Tung},
       title={On the regularity of the {M}onge-{A}mp\`ere equation
  {$\det(\partial^2u/\partial x_i\partial x_j)=F(x,u)$}},
        date={1977},
     journal={Comm. Pure Appl. Math.},
      volume={30},
      number={1},
       pages={41\ndash 68},
}

\bib{cohn-vossen1927}{article}{
      author={Cohn-Vossen, Stefan},
       title={Zwei {S}\"atze \"uber die {S}tarrheit der {E}ifl\"achen},
        date={1927},
     journal={Nachr. Ges. Wiss. G\"ottingen, Math.-Phys. Kl.},
       pages={125\ndash 134},
}

\bib{cohn-vossen1930}{article}{
      author={Cohn-Vossen, Stefan},
       title={Unstarre geschlossene {F}l\"achen},
        date={1930},
     journal={Math. Ann.},
      volume={102},
       pages={10\ndash 29},
}

\bib{connelly1977}{article}{
      author={Connelly, Robert},
       title={A counterexample to the rigidity conjecture for polyhedra},
        date={1977},
     journal={Publ. Math. Inst. Hautes \'Etudes Sci.},
      volume={47},
       pages={333\ndash 338},
}

\bib{ghomi2002}{article}{
      author={Ghomi, Mohammad},
       title={Gauss map, topology, and convexity of hypersurfaces with
  nonvanishing curvature},
        date={2002},
        ISSN={0040-9383},
     journal={Topology},
      volume={41},
      number={1},
       pages={107\ndash 117},
      review={\MR{1871243}},
}

\bib{ghomi-qiu-raffaelli-zhang2026}{misc}{
      author={Ghomi, Mohammad},
      author={Qiu, Hongda},
      author={Raffaelli, Matteo},
      author={Zhang, Yuetong},
       title={{TightSurfaces.nb}},
        date={2026},
        note={Mathematica notebook, available at
  \url{https://ghomi.math.gatech.edu/MathematicaNBs/TightSurfaces.nb}},
}

\bib{ghomi-raffaelli2025}{article}{
      author={Ghomi, Mohammad},
      author={Raffaelli, Matteo},
       title={Topology of closed asymptotic curves on negatively curved
  surfaces},
        date={2025},
        ISSN={1050-6926,1559-002X},
     journal={J. Geom. Anal.},
      volume={35},
      number={12},
       pages={Paper No. 381, 12 pp.},
         url={https://doi.org/10.1007/s12220-025-02200-3},
      review={\MR{4970320}},
}

\bib{ghomi-spruck2020}{article}{
      author={Ghomi, Mohammad},
      author={Spruck, Joel},
       title={Rigidity of nonnegatively curved surfaces relative to a curve},
        date={2020},
        ISSN={1073-7928,1687-0247},
     journal={Int. Math. Res. Not. IMRN},
      number={17},
       pages={5387\ndash 5400},
         url={https://doi.org/10.1093/imrn/rny167},
      review={\MR{4146341}},
}

\bib{han-hong2006}{book}{
      author={Han, Qing},
      author={Hong, Jia-Xing},
       title={Isometric embedding of {R}iemannian manifolds in {E}uclidean
  spaces},
      series={Mathematical Surveys and Monographs},
   publisher={American Mathematical Society},
     address={Providence, RI},
        date={2006},
      volume={130},
        ISBN={978-0-8218-4071-9; 0-8218-4071-1},
      review={\MR{2261749}},
}

\bib{han-hong-huang2017}{article}{
      author={Han, Qing},
      author={Hong, Jiaxing},
      author={Huang, Genggeng},
       title={Compactness of {A}lexandrov-{N}irenberg surfaces},
        date={2017},
     journal={Comm. Pure Appl. Math.},
      volume={70},
      number={9},
       pages={1706\ndash 1753},
      review={\MR{3684308}},
}

\bib{han-khuri2011}{article}{
      author={Han, Qing},
      author={Khuri, Marcus},
       title={Rigidity in the class of orientable compact surfaces of minimal
  total absolute curvature},
        date={2011},
        ISSN={0926-2245,1872-6984},
     journal={Differential Geom. Appl.},
      volume={29},
      number={4},
       pages={463\ndash 472},
         url={https://doi.org/10.1016/j.difgeo.2011.04.035},
      review={\MR{2811658}},
}

\bib{kovaleva1968}{article}{
      author={Kovaleva, G.~A.},
       title={Example of a surface that is homotopic to a tube and has a closed
  asymptotic line},
        date={1968},
     journal={Math. Notes},
      volume={3},
      number={4},
       pages={257\ndash 263},
         url={https://doi.org/10.1007/BF01159941},
}

\bib{lelieuvre1888}{article}{
      author={Lelieuvre, M.},
       title={Sur les lignes asymptotiques et leur repr\'esentation
  sph\'erique},
        date={1888},
     journal={Bull. Sci. Math. (2)},
      volume={12},
       pages={126\ndash 128},
         url={https://gallica.bnf.fr/ark:/12148/bpt6k486248d/f211},
}

\bib{ni-zhang-zhang2026}{misc}{
      author={Ni, Lei},
      author={Zhang, Wei},
      author={Zhang, Yijian},
       title={On questions of {P}ogorelov and {T}oponogov},
        date={2026},
         url={https://doi.org/10.48550/arXiv.2606.29231},
        note={arXiv:2606.29231 [math.DG]},
}

\bib{nirenberg1963}{incollection}{
      author={Nirenberg, Louis},
       title={Rigidity of a class of closed surfaces},
        date={1963},
   booktitle={Nonlinear {P}roblems ({P}roc. {S}ympos., {M}adison, {W}is.,
  1962)},
   publisher={Univ. of Wisconsin Press, Madison, Wis.},
       pages={177\ndash 193},
      review={\MR{0150705}},
}

\bib{pogorelov1973}{book}{
      author={Pogorelov, A.~V.},
       title={Extrinsic geometry of convex surfaces},
   publisher={American Mathematical Society},
     address={Providence, R.I.},
        date={1973},
        note={Translated from the Russian by Israel Program for Scientific
  Translations, Translations of Mathematical Monographs, Vol. 35},
      review={\MR{0346714}},
}

\bib{pogorelov1978}{book}{
      author={Pogorelov, A.~V.},
       title={The {M}inkowski multidimensional problem},
   publisher={V. H. Winston \& Sons, Washington, D.C.},
        date={1978},
        note={Translated from the Russian by Vladimir Oliker, Scripta Series in
  Mathematics},
}

\bib{rembs1952}{article}{
      author={Rembs, Eduard},
       title={Zur {V}erbiegung von {F}l\"achen im {G}rossen},
        date={1952},
        ISSN={0025-5874},
     journal={Math. Z.},
      volume={56},
       pages={271\ndash 279},
      review={\MR{0054308}},
}

\bib{spivak:v5}{book}{
      author={Spivak, Michael},
       title={A comprehensive introduction to differential geometry. {V}ol.
  {V}},
     edition={Second},
   publisher={Publish or Perish Inc.},
     address={Wilmington, Del.},
        date={1979},
        ISBN={0-914098-83-7},
      review={\MR{532834}},
}

\bib{trudinger-wang2008}{incollection}{
      author={Trudinger, Neil~S.},
      author={Wang, Xu-Jia},
       title={The {M}onge-{A}mp\`ere equation and its geometric applications},
        date={2008},
   booktitle={Handbook of geometric analysis. {N}o. 1},
      series={Adv. Lect. Math. (ALM)},
      volume={7},
   publisher={Int. Press, Somerville, MA},
       pages={467\ndash 524},
}

\bib{yau2000}{incollection}{
      author={Yau, S.-T.},
       title={Review of geometry and analysis},
        date={2000},
   booktitle={Mathematics: frontiers and perspectives},
   publisher={Amer. Math. Soc.},
     address={Providence, RI},
       pages={353\ndash 401},
      review={\MR{1754787}},
}

\end{biblist}
\end{bibdiv}

\end{document}